\documentclass[letterpaper,11pt]{article}

\usepackage{amssymb,amsmath,amscd,amsfonts}
\usepackage{epsf,epsfig}
\usepackage{times,bm,mathrsfs,bbm}
\usepackage{color,graphicx}

\definecolor{Red}{rgb}{1,0,0}
\definecolor{Blue}{rgb}{0,0,1}
\definecolor{Olive}{rgb}{0.41,0.55,0.13}
\definecolor{Green}{rgb}{0,1,0}
\definecolor{MGreen}{rgb}{0,0.8,0}
\definecolor{DGreen}{rgb}{0,0.55,0}
\definecolor{Yellow}{rgb}{1,1,0}
\definecolor{Cyan}{rgb}{0,1,1}
\definecolor{Magenta}{rgb}{1,0,1}
\definecolor{Orange}{rgb}{1,.5,0}
\definecolor{Violet}{rgb}{.5,0,.5}
\definecolor{Purple}{rgb}{.75,0,.25}
\definecolor{Brown}{rgb}{.75,.5,.25}
\definecolor{Grey}{rgb}{.5,.5,.5}
\definecolor{Black}{rgb}{0,0,0}

\newtheorem{theorem}{Theorem}
\newtheorem{proposition}{Proposition}

\newtheorem{definition}{Definition}
\newtheorem{example}{Example}
\newtheorem{lemma}{Lemma}

\newtheorem{remark}{Remark}

\newenvironment{proof}{\noindent{\textbf{Proof:}}}{$\blacksquare$\vskip\belowdisplayskip}

\newcommand{\acal}{\mathcal{A}}
\newcommand{\bcal}{\mathcal{B}}

\newcommand{\ecal}{\mathcal{E}}
\newcommand{\fcal}{\mathcal{F}}

\newcommand{\lcal}{\mathcal{L}}

\newcommand{\qcal}{\mathcal{Q}}
\newcommand{\rcal}{\mathcal{R}}

\newcommand{\tcal}{\mathcal{T}}

\newcommand{\zcal}{\mathcal{Z}}

\newcommand{\ind}{\mathbbm{1}}

\newcommand{\eps}{\varepsilon}

\newcommand{\prob}{\mathbb{P}}

\newcommand{\poly}{\mathrm{poly}}

\newcommand{\weight}{\tau}
\newcommand{\eweight}{\hat\tau}
\newcommand{\path}{\mathrm{Path}}
\newcommand{\dist}{\mathrm{d}}

\newcommand{\phy}{\tcal}
\newcommand{\rt}{\rho}
\newcommand{\hmgt}[1]{T^{(#1)}}
\newcommand{\hmgv}[1]{V^{(#1)}}
\newcommand{\hmgl}[1]{L^{(#1)}}
\newcommand{\hmgll}[2]{L^{(#1)}_{#2}}
\newcommand{\hmge}[1]{E^{(#1)}}

\newcommand{\hmgphy}[1]{\phy^{(#1)}}
\newcommand{\hmgrt}[1]{\rt^{(#1)}}
\newcommand{\states}{[q]}
\newcommand{\nstates}{q}
\newcommand{\nstatespotts}{q}
\newcommand{\qq}{\nstatespotts}
\newcommand{\QQ}{R}
\newcommand{\trees}{\mathbb{T}}

\newcommand{\rates}{\mathbb{Q}}
\newcommand{\gtr}{\mathbb{G}}
\newcommand{\law}{\lcal}
\newcommand{\sphy}{\mathbb{Y}}
\newcommand{\shmgphy}{\mathbb{HY}}

\newcommand{\deep}{\overline{\mathbb{FP}}}

\newcommand{\overm}[2]{\overline{M}^{#1,#2}}

\newcommand{\critks}{g^{\star}_{\mathrm{Lin}}}
\newcommand{\critksq}{g_{\mathrm{Lin}}(Q)}
\newcommand{\critmlq}{g_{\mathrm{ML}}(Q)}
\newcommand{\critperc}{g^{\star}_{\mathrm{Perc}}}

\newcommand{\critlargeq}{g^+_q}
\newcommand{\weightmax}{\weight^+}

\newcommand{\bias}{\mathrm{b}}
\newcommand{\maxbias}{\bar{\bias}}
\newcommand{\weightb}{\weight_\bias}

\newenvironment{app-proof}[1]{\noindent{\textbf{Proof of #1:}}}{$\blacksquare$\vskip\belowdisplayskip}

\begin{document}

\title{\vspace{-3cm}
On the inference of large phylogenies with long branches:
How long is too long?
}
\author{
Elchanan Mossel\thanks{Weizmann Institute and U.C. Berkeley. 
Supported by NSF Career Award (DMS 054829),
by ONR award N00014-07-1-0506, by ISF grant 1300/08 and by Marie
Curie grant   PIRG04-GA-2008-239317}\and
S\'ebastien Roch\thanks{UCLA.}\and
Allan Sly\thanks{Microsoft Research.}
}
\maketitle

\begin{abstract}
The accurate reconstruction of phylogenies from short molecular sequen\-ces
is an important problem in computational biology. Recent work has highlighted
deep connections between sequence-length requirements
for high-probability phylogeny reconstruction
and the related problem of the estimation of ancestral sequences.
In [Daskalakis et al.'09], building on the work of [Mossel'04], a tight
sequence-length requirement was obtained for the simple CFN model of substitution, that is,
the case of a two-state symmetric rate matrix $Q$.
In particular the required sequence length for high-probability reconstruction was
shown to undergo a sharp transition (from $O(\log n)$ to
$\hbox{poly}(n)$, where $n$ is the number of leaves)
at the ``critical'' branch length $\critmlq$ (if it exists) of the
ancestral reconstruction problem defined roughly as follows:
below $\critmlq$ the sequence at the
root can be accurately estimated from sequences at the leaves on deep
trees, whereas above $\critmlq$ information decays exponentially quickly
down the tree.

Here we consider a more general evolutionary model, the
GTR model, where the $\nstates\times\nstates$ rate matrix $Q$ is reversible
with $\nstates \geq 2$. For this model, recent results of [Roch'09]
show that the tree can be accurately reconstructed
with sequences of length $O(\log(n))$ when the branch lengths are below
$\critksq$, known as the Kesten-Stigum (KS) bound, up to which
ancestral sequences can be accurately estimated using simple linear
estimators.
Although for the CFN model $\critmlq = \critksq$  (in other words, linear
ancestral estimators are in some sense best possible),
it is known that for the more general GTR models one has $\critmlq \geq \critksq$
with a \emph{strict} inequality in many cases.
Here, we show that this phenomenon also holds for phylogenetic reconstruction by
exhibiting a family of symmetric models $Q$ and a phylogenetic reconstruction
algorithm which recovers the tree from $O(\log n)$-length sequences for some branch lengths
in the range $(\critksq,\critmlq)$. Second we prove that phylogenetic
reconstruction under GTR models requires a polynomial sequence-length for branch lengths
above $\critmlq$.

\end{abstract}

\section{Introduction}

\paragraph{Background.}
Recent years have witnessed a convergence of models and problems
from evolutionary biology, statistical physics, and computer science.
Standard stochastic models of molecular evolution, such as the
Cavender-Farris-Neyman (CFN) model (a.k.a.~the Ising model or Binary Symmetric Channel (BSC))
or the Jukes-Cantor (JC) model (a.k.a.~the Potts model),
have been extensively studied from all these different perspectives
and fruitful insights have emerged, notably in the area of
computational phylogenetics.

Phylogenetics~\cite{SempleSteel:03,Felsenstein:04} is centered around the reconstruction of evolutionary histories
from molecular data extracted from modern species.
The assumption is that molecular data consists of aligned sequences and
that each position in the sequences evolves independently according to a
Markov model on a tree, where the key parameters are
(see Section~\ref{section:definitions} for formal definitions):
\begin{itemize}
\item {\em Rate matrix.}
A $\nstates \times \nstates$ mutation rate matrix $Q$, where $\nstates$ is the alphabet size.
A typical alphabet is the set of nucleotides $\{\mathrm{A},\mathrm{C},\mathrm{G},\mathrm{T}\}$,
but here we allow more general state spaces.
Without loss of generality, we denote the alphabet by $\states = \{1,\ldots,\nstates\}$.
The $(i,j)$'th entry of $Q$ encodes the rate at which state $i$ mutates into state $j$.
\item {\em Tree.}
An evolutionary tree $T$, where the leaves are the modern species and each branching represents a past speciation event.
We denote the leaves by $[n] = \{1,\ldots,n\}$.
\item {\em Branch lengths.}
For each edge $e$,
we have a scalar branch length $\weight(e)$
which measures the expected total number of substitutions per site along edge $e$.
Roughly speaking, $\weight(e)$ is the time duration between the end points of $e$
multiplied by the mutation rate.
\end{itemize}
We consider the following two closely related problems:
\begin{enumerate}
\item
{\bf Phylogenetic Tree Reconstruction (PTR).}
Given $n$ molecular sequences of length $k$ (one for each leaf)
\begin{equation*}
\{s_a = (s^i_a)_{i=1}^k\}_{a\in[n]}
\end{equation*}
with $s^i_a \in [q]$,
which have evolved according to the process above with independent sites,
reconstruct the topology of the evolutionary tree.
\item
{\bf Ancestral State Reconstruction (ASR).}
Given a fully specified rooted tree and a single state $s^1_a$ at each leaf $a$ of the tree, estimate (better than ``random") the state at the root of the tree, independently of the depth of the tree.
\end{enumerate}
In both cases, longer edge lengths correspond to more mutations---and hence more noise---making both reconstruction problems more challenging. Our overriding goal is to extend efficient phylogenetic reconstruction to trees with as large branch lengths as possible.

\paragraph{Reconstruction thresholds.}
Alternatively, the second problem can be interpreted
in terms of correlation decay along the tree
or as a broadcasting problem on a tree-network.
It has thus been extensively studied in statistical physics, probability theory, and computer science.
See e.g.~\cite{EvKePeSc:00} and references therein.
A crucial parameter in the ASR problem is $\weightmax(T) = \max_e \weight(e)$,
the maximal branch length in the tree.

One class of ancestral estimators is particularly well understood, the so-called linear estimators.
See Section~\ref{section:definitions} for a formal definition. In essence, linear estimators
are simply a form of weighted majority. In~\cite{MosselPeres:03},
it was shown that there exists a critical parameter
$\critksq = \lambda_Q^{-1} \ln\sqrt{2}$, where $-\lambda_Q$ is the largest negative eigenvalue
of the rate matrix $Q$, such that:
\begin{itemize}
\item if $\weightmax < \critksq$, for all trees with $\weightmax(T) = \weightmax$
a \emph{well-chosen} linear estimator provides a good solution to the ASR,
\item if $\weightmax > \critksq$, there exist trees with $\weightmax(T) = \weightmax$
for which ASR is impossible for \emph{any} linear estimator,
that is, the correlation between the best linear root estimate and
the true root value decays exponentially in the depth of the tree.
\end{itemize}
For formal definitions, see~\cite{MosselPeres:03}. The threshold $\critksq = \lambda_Q^{-1} \ln\sqrt{2}$ is also known to
be the critical threshold for {\em robust (ancestral) reconstruction}, see~\cite{JansonMossel:04} for details.

For more general ancestral estimators, only partial results are known.
For the two-state symmetric $Q$ (the CFN model), impossibility
of reconstruction as above holds, when $\weightmax(T) > \critksq$, not only
for linear estimators but also for {\em any} estimator, including for
instance maximum likelihood. In other words, for the CFN model
linear estimators are in some sense best possible. This phenomenon
also holds for symmetric models (i.e., where all non-diagonal entries of $Q$
are identical) with $q=3$ states~\cite{Sly:09} (at least, for high degree trees).
However, for symmetric models on $q\geq 5$ states,
it is known that ASR is possible beyond $\critksq$,
up to a critical branch length $\critmlq > \critksq$
which is not known explicitly~\cite{Mossel:01,Sly:09}.
Larger values of $q$ here correspond for instance to models of protein evolution.
ASR beyond $\critksq$ can be achieved with a maximum likelihood estimator although
in some cases special estimators have been devised (for instance, 
symmetric models with large $q$)~\cite{Mossel:01}.
In this context, $\critksq$ is refered to as the
{\em Kesten-Stigum bound}~\cite{KestenStigum:67}.
We sometimes call the condition $\weightmax(T) < \critksq$ the ``KS phase''
and the condition $\weightmax(T) < \critmlq$ the ``reconstruction phase.''
%

For general reversible rate matrices, it is not even known whether there is a {\em unique}
reconstruction threshold $\critmlq$ such that
ASR is possible for $\weightmax(T) < \critmlq$ and impossible for $\weightmax(T) > \critmlq$.
The general question of finding the threshold $\critmlq$ for ASR is extremely challenging and has been answered for only a very small number of channels.


\paragraph{Steel's Conjecture.}
A striking conjecture of Steel~\cite{Steel:01}
postulates a deep connection between PTR and ASR.
More specifically, the conjecture states that for CFN models if $\weightmax(T) < \critksq$
then PTR can be achieved with sequence length $k = O(\log n)$.  
This says that, when we can accurately estimate the states 
of vertices deep inside a \emph{known} tree, then it 
is also possible to accurately reconstruct the topology of an \emph{unknown} tree
with very short sequence lengths.

In fact, since the number of trees on $n$ labelled leaves is $2^{\Theta(n \log n)}$,
this is an optimal sequence length up to constant factors---that is, we cannot
hope to distinguish so many trees with fewer potential datasets.
The proof of Steel's conjecture was established in~\cite{Mossel:04a}
for balanced trees and in~\cite{DaMoRo:09} for general (under the additional
assumption that branch lengths are discretized).
Furthermore, results of Mossel~\cite{Mossel:03,Mossel:04a}
show that for $\weightmax(T) > \critksq$ a polynomial sequence length is needed
for correct phylogenetic reconstruction.
For symmetric models, the results of~\cite{Mossel:04a,DaMoRo:09} imply that
it is possible to reconstruct phylogenetic trees from sequences of length
$O(\log n)$ when $\weightmax(T) < \critksq$.
These results cover classical models such as the JC model ($\nstates = 4$).
Recent results of Roch~\cite{Roch:09}, building on~\cite{Roch:08,PeresRoch:09},
show that for any reversible mutation matrix $Q$,
it is possible to reconstruct phylogenetic trees from
$O(\log(n))$-length sequences again when $\weightmax(T) < \critksq$.

However, these results leave the following important problem open:
\begin{itemize}
\item
As we mentioned before, for symmetric models on $q\geq 5$ states,
it is known that ASR is possible for $\weightmax(T) < \critmlq$,
where $\critmlq > \critksq$.
A natural question is to ask if the ``threshold'' for PTR is $\critmlq$ (i.e., the threshold for ASR)
or $\critksq$ or perhaps another value.
(Note that for the CFN model, the threshold for PTR has been shown to be
$\critksq$ but in that case it so happens that $\critksq = \critmlq$.)
\end{itemize}

\paragraph{Our contributions.}
Our main results are the following:
\begin{itemize}
\item
We show that for symmetric models $Q$ with large $\nstates$, it is possible to reconstruct phylogenetic trees with
$O(\log n)$-length sequences whenever
$\weightmax(T) < \critlargeq$ where $\critksq < \critlargeq < \critmlq$.
We thus show that PTR from logarithmic sequences is sometimes possible for branch lengths {\em above} the KS bound.
\item
We also show how to generalize the arguments of~\cite{Mossel:03,Mossel:04a}
to show that for any $Q$ and $\weightmax(T) > \critmlq$
it holds that correct phylogenetic reconstruction requires polynomial-length sequences in general.
The same idea is used in~\cite{Mossel:03,Mossel:04a} and the argument presented here. The main difference is that in the arguments
in~\cite{Mossel:03,Mossel:04a} used mutual information together with coupling while the more elegant argument presented here uses coupling only. The results of~\cite{Mossel:03} apply for general models but are not tight even for the CFN model.
The argument in~\cite{Mossel:04a} gives tight results for the CFN model. It is possible to extend that argument to more general models, but we prefer the simpler proof given in the current paper.
\end{itemize}


\paragraph{Organization.}
We begin with preliminaries and the formal statements of our results in Section~\ref{section:definitions}.
The proof of our upper bound can be found in Section~\ref{section:upper}.
The proof of our lower bound can be found in Section~\ref{section:lower}.

\section{Definitions and Results}\label{section:definitions}

\subsection{Basic Definitions}

\noindent\textbf{Phylogenies.} We define phylogenies and evolutionary distances more formally.
\begin{definition}[Phylogeny]
A {\em phylogeny} is a rooted, edge-weighted, leaf-labeled tree
$\phy = (V,E,[n],\rt;\weight)$
where:
$V$ is the set of vertices;
$E$ is the set of edges;
$L = [n] = \{1,\ldots,n\}$ is the set of leaves;
$\rt$ is the root;
$\weight : E \to (0,+\infty)$ is a positive edge weight function.
We further assume that all internal nodes in $\phy$ have degree $3$
except for the root $\rt$ which has degree $2$.
We let $\sphy_n$ be the set of all such phylogenies on $n$ leaves
and we denote $\sphy = \{\sphy_n\}_{n\geq 1}$.
\end{definition}
\begin{definition}[Tree Metric]
For two leaves $a,b \in [n]$, we denote by $\path(a,b)$
the set of edges on the unique path between $a$ and $b$.
A {\em tree metric} on a set $[n]$ is a positive function
$\dist:[n]\times[n] \to (0,+\infty)$ such that there exists
a tree $T = (V,E)$ with leaf set $[n]$ and an edge weight
function $w:E \to (0,+\infty)$ satisfying the following: for all leaves
$a,b \in [n]$
\begin{equation*}
\dist(a,b) = \sum_{e\in \path(a,b)} w_e.
\end{equation*}
For convenience, we denote by $\left(\weight(a,b)\right)_{a,b\in [n]}$
the tree metric corresponding to phylogeny $\phy = (V,E,[n],\rt;\weight)$.
We extend $\weight(u,v)$ to all vertices $u,v \in V$ in the obvious way.
\end{definition}
\begin{example}[Homogeneous Tree]\label{ex:homo}
For an integer $h \geq 0$, we denote by
$\hmgphy{h} = (\hmgv{h}, \hmge{h}, \hmgl{h}, \hmgrt{h}; \weight)$
a rooted phylogeny where $\hmgt{h}$ is the $h$-level
complete binary tree with arbitrary edge weight function
$\weight$ and $\hmgl{h} = [2^h]$.
For $0\leq h'\leq h$, we let $\hmgll{h}{h'}$ be the
vertices on level $h - h'$ (from the root). In particular,
$\hmgll{h}{0} = \hmgl{h}$ and $\hmgll{h}{h} = \{\hmgrt{h}\}$.
We let $\shmgphy = \{\shmgphy_n\}_{n\geq 1}$ be the set of all phylogenies
with homogeneous underlying trees.
\end{example}

\noindent\textbf{Model of molecular sequence evolution.}
Phylogenies are reconstructed from molecular sequences
extracted from the observed species. The standard model
of evolution for such sequences is a Markov model
on a tree (MMT).
\begin{definition}[Markov Model on a Tree]
Let $\nstates \geq 2$. Let $n \geq 1$ and let $T = (V,E,[n],\rt)$
be a rooted tree with leaves labeled in $[n]$.
For each edge $e \in E$, we are given a $\nstates\times\nstates$
stochastic matrix
$M^e = (M^e_{ij})_{i,j \in \states}$,
with fixed stationary distribution
$\pi = (\pi_i)_{i\in \states}$.
An MMT $(\{M^e\}_{e\in E}, T)$
associates a state $s_v$ in $\states$ to each
vertex $v$ in $V$ as follows:
pick a state for the root $\rt$ according to $\pi$;
moving away from the root, choose a state for each
vertex $v$ independently according to the distribution
$(M^e_{s_u, j})_{j\in\states}$,
with $e = (u,v)$ where $u$ is the parent of $v$.
\end{definition}
The most common MMT used in phylogenetics is the
so-called general time-reversible (GTR) model.
\begin{definition}[GTR Model]
Let $\states$ be a set of character states with $\nstates = |\states|$
and $\pi$ be a distribution on $\states$ satisfying $\pi_i > 0$ for all
$i\in\states$.
For $n \geq 1$, let $\phy = (V,E,[n],\rt;\weight)$ be a phylogeny.
Let $Q$ be a $\nstates\times\nstates$ rate matrix, that is,
$Q_{ij} > 0$ for all $i\neq j$ and
$\sum_{j\in \states} Q_{ij} = 0$,
for all $i \in \states$.
Assume $Q$ is reversible with respect to $\pi$, that is,
$\pi_i Q_{ij} = \pi_j Q_{ji}$,
for all $i,j \in \states$.
The GTR model on $\phy$ with rate matrix $Q$ is an MMT
on $T = (V,E,[n], \rt)$ with transition matrices
$M^e = e^{\weight_e Q}$,
for all $e\in E$.
By the reversibility assumption, $Q$ has $\nstates$
real eigenvalues
$0 = \Lambda_1 > \Lambda_2 \geq \cdots \geq \Lambda_{\nstates}$.
We normalize $Q$ by fixing $\Lambda_2 = -1$.
We denote by $\rates_\nstates$ the set of all such rate matrices.
We let $\gtr_{n,\nstates} = \sphy_n \otimes \rates_\nstates$
be the set of all $\nstates$-state GTR models on
$n$ leaves. We denote $\gtr_\nstates =
\left\{\gtr_{n,\nstates}\right\}_{n \geq 1}$.
We denote by $s_W$ the vector of states on the vertices
$W\subseteq V$. In particular, $s_{[n]}$ are the states
at the leaves.
We denote by $\law_{\phy,Q}$ the distribution of $s_{[n]}$.
\end{definition}

GTR models are often used in their full generality in the biology literature,
but they also encompass several popular special cases such as the CFN model and the JC model.
\begin{example}[$\qq$-State Symmetric Model]\label{ex:symmetric}
The {\em $\qq$-state Symmetric model} (also called $\qq$-state Potts model) is the GTR model with $\qq \geq 2$ states,
$\pi = (1/\qq,\ldots, 1/\qq)$,
and $Q = Q^{(\qq)}$ where
\begin{equation*}
Q^{(\qq)}_{ij}
=
\left\{
\begin{array}{ll}
-\frac{\qq - 1}{\qq} & \mbox{if $i=j$}\\
\frac{1}{\qq} & \mbox{o.w.}
\end{array}
\right.
\end{equation*}
It is easy to check that $\Lambda_2(Q) = -1$.
The special cases $\qq=2$ and $\qq=4$ are called respectively the CFN and JC models in the biology literature.
We denote their rate matrices by $Q^{\mathrm{CFN}}, Q^{\mathrm{JC}}$.
For an edge $e$ of length $\weight_e > 0$, let
\begin{equation*}
\delta_e = \frac{1}{\qq}\left(1 - e^{-\weight_e}\right).
\end{equation*}
Then, we have
\begin{equation*}
(M_e)_{ij} = (e^{\weight_e Q})_{ij}
=
\left\{
\begin{array}{ll}
1 -(\qq - 1) \delta_e & \mbox{if $i=j$}\\
\delta_e & \mbox{o.w.}
\end{array}
\right.
\end{equation*}
\end{example}

\noindent\textbf{Phylogenetic reconstruction.}
A standard assumption in molecular evolution
is that each site in a sequence (DNA, protein, etc.)
evolves {\em independently} according to a
Markov model on a tree, such as the GTR model above.
Because of the reversibility assumption, the root
of the phylogeny cannot be identified and we reconstruct phylogenies
up to their root.
\begin{definition}[Phylogenetic Reconstruction Problem]
Let $\widetilde\sphy = \{\widetilde\sphy_n\}_{n\geq 1}$
be a subset of phylogenies and
$\widetilde\rates_\nstates$ be a subset of rate matrices on $\nstates$ states.
Let $\phy = (V,E,[n],\rt;\weight) \in \widetilde\sphy$.
If $T = (V,E,[n],\rt)$ is the rooted tree underlying $\phy$,
we denote by $T_{-}[\phy]$ the tree $T$ where the root
is removed: that is, we replace the two edges adjacent to the root
by a single edge. We denote by $\trees_n$ the set of all
leaf-labeled trees on $n$ leaves with internal degrees $3$
and we let $\trees = \{\trees_n\}_{n\geq 1}$.
A {\em phylogenetic
reconstruction algorithm} is a collection of maps
$\acal = \{\acal_{n,k}\}_{n,k \geq 1}$
from sequences $(s^i_{[n]})_{i=1}^k \in (\states^{[n]})^k$
to leaf-labeled trees $T \in \trees_n$.
We only consider algorithms $\acal$ computable in time polynomial
in $n$ and $k$.
Let $k(n)$ be an increasing function of $n$. We say that $\acal$
solves the {\em phylogenetic reconstruction problem}
on $\widetilde\sphy \otimes \widetilde\rates_\nstates$ with sequence length $k = k(n)$
if for all
$\delta > 0$, there is $n_0 \geq 1$ such that for all $n \geq n_0$,
$\phy \in \widetilde\sphy_n$, $Q \in \widetilde\rates_\nstates$,
\begin{equation*}
\prob\left[\acal_{n,k(n)}\left((s^i_{[n]})_{i=1}^{k(n)}\right) =
T_-[\phy]\right] \geq 1 - \delta,
\end{equation*}
where $(s^i_{[n]})_{i=1}^{k(n)}$ are i.i.d.~samples from $\law_{\phy,Q}$.
\end{definition}
An important result of this kind was given by Erdos et al.~\cite{ErStSzWa:99a}.
Let $\alpha\geq 1$ and $\nstates \geq 2$. The set of rate matrices $Q \in \rates_\nstates$
such that $\mathrm{tr}(Q) \geq -\alpha$ is denoted $\rates_{\nstates,\alpha}$.
Let $0 < f < g < +\infty$ and denote by $\sphy^{f,g}$ the set of
all phylogenies $\phy = (V,E,[n],\rt;\weight)$ satisfying
$f < \weight_e < g,\ \forall e\in E$.
Then, Erdos et al.~showed (as rephrased in our setup) that,
for all $\alpha \geq q-1$, $\nstates \geq 2$, and all $0 < f < g < +\infty$,
the phylogenetic reconstruction problem on $\sphy^{f,g}\otimes\rates_{\nstates,\alpha}$
can be solved with $k = \poly(n)$.
(In fact, they proved a more general result allowing rate matrices
to vary across different edges.)
In the case of the Potts model, this result was improved by
Daskalakis et al.~\cite{DaMoRo:09} (building on~\cite{Mossel:04a})
in the Kesten-Stigum (KS) reconstruction phase, that is, when
$g < \critksq = \critks \equiv \ln\sqrt{2}$.
They showed that, for all $0 < f < g < \critks$,
the phylogenetic reconstruction problem on $\sphy^{f,g}\otimes \{Q^{(\qq)}\}$
can be solved with $k = O(\log(n))$. More recently, the latter result was extended 
to GTR models by Roch~\cite{Roch:09}, building on~\cite{Roch:08,PeresRoch:09}. 
But prior to our work, no PTR algorithm had been shown to extend beyond  $\critks$.

\subsection{Our Results}

\noindent\textbf{Positive result.}
In our first result, we extend
logarithmic reconstruction results for $\qq$-state symmetric models
to $\ln\sqrt{2} < g < \ln 2$ for large enough $\qq$.
This is the first result of this type beyond the KS bound.
\begin{theorem}[Logarithmic Reconstruction beyond the KS Transition]\label{thm:potts}
Let $0 < f < g < +\infty$ and denote by $\shmgphy^{f,g}$ the set of
all homogeneous phylogenies $\phy = (V,E,[n],\rt;\weight)$ satisfying
$f < \weight_e < g,\ \forall e\in E$.
Let $\critperc = \ln 2$.
Then, for all $0 < f < g < \critperc$,
there is $\QQ \geq 2$ such that for all $\qq > \QQ$
the phylogenetic reconstruction problem on $\shmgphy^{f,g}\otimes\{Q^{(\qq)}\}$
can be solved with $k = O(\log(n))$.
\end{theorem}
Theorem~\ref{thm:potts} can be extended to general phylogenies using the techniques
of~\cite{DaMoRo:09}, although then one requires discretized branch lengths.
See~\cite{DaMoRo:09} for details.

\paragraph{Negative result.}
In our second result, we show that for $g > \critmlq$ the number of samples $k$ must grow
polynomially in $n$. In particular, this is true for the $\qq$-state symmetric model
for all $\qq \geq 2$ and $g > \ln 2$ by the results of~\cite{Mossel:01}.
\begin{theorem}[Polynomial Lower Bound Above $\critmlq$ (see also~\cite{Mossel:03,Mossel:04a})]\label{thm:lower}
Let $Q\in \rates_\nstates$ and $f = g > \critmlq$. Then the phylogenetic reconstruction problem
on $\shmgphy^{f,g}\otimes\{Q\}$ requires $k = \Omega(n^{\alpha})$ for some $\alpha > 0$
(even assuming $Q$ and $g$ are known exactly beforehand).
\end{theorem}

\begin{remark}[Biological Convention]
Our normalization of $Q$ differs from standard biological
convention where it is assumed that the total rate of change
per unit time at stationarity is 1, that is,
\begin{equation*}
\sum_{i} \pi_i Q_{ii} = -1.
\end{equation*}
See e.g.~\cite{Felsenstein:04}.
Let $-\lambda_Q$ denote the largest negative eigenvalue under this
convention. Then, the Kesten-Stigum bound is given by the solution to
\begin{equation*}
2 e^{-2\lambda_Q \critksq} = 1.
\end{equation*}
For instance, in the Jukes-Cantor model one has
\begin{equation*}
\critksq = \frac{3}{8}\ln 2.
\end{equation*}
\end{remark}

\section{Upper Bound for Large $\qq$}\label{section:upper}

\subsection{Root Estimator}

The basic ingredient behind logarithmic reconstruction results is an accurate
estimator of the root state. In the KS phase, this can be achieved by
majority-type procedures. See~\cite{Mossel:98,EvKePeSc:00,Mossel:04a}.
In the reconstruction phase beyond the KS phase however, a more sophisticated
estimator is needed.  
In this subsection we define an accurate root estimator which 
does not depend on the edge lengths.

\paragraph{Random Cluster Methods.}
We use a convenient percolation representation of the
ferromagnetic Potts model on trees. Let $\qq \geq 2$
and $\phy = (V,E,[n],\rt;\weight) \in \shmgphy_n$
with corresponding $(\delta_e)_{e\in E}$.
Run a percolation process on $T = (V,E)$
where edge $e$ is open with probability
$1 - \qq \delta_e$. Then associate
to each open cluster a state according
to the uniform distribution on $\states$.
The state so obtained $(s_v)_{v\in V}$ has the same distribution
as the GTR model $(\phy,Q^{(\qq)})$.

We will use the following definition.
Let $T'$ be a subtree of $T$ which is rooted at
$\rt$. We say that $T'$ is an \emph{$l$-diluted binary tree}
if, for all $s$, all the vertices of $T'$ at level $sl$
have exactly $2$ descendants at level $(s+1)l$.
(Assume for now that $\log_2 n$ is a multiple of $l$.)
For a state $i\in \states$ and assignment $s_{[n]}$ at the leaves,
we say that the event $\bcal_{i,l}$ holds if there is
a $l$-diluted binary tree with state $i$ at all its leaves according
to $s_{[n]}$.
Let $B_l$ be the set of all $i$ such that
$\bcal_{i,l}$ holds. Consider the following estimator:
pick a state $X$ uniformly at random in $\states$ and let
\begin{equation*}
\bar{s}^l_\rt
=
\left\{
\begin{array}{ll}
X, & \mbox{if $X \in B_l$}\\
\mbox{pick uniformly in $\states - \{X\}$,} & \mbox{o.w.}
\end{array}
\right.
\end{equation*}
We use the following convention. If $\log_2 n$ is not a multiple
of $l$, we add levels of $0$-length edges to $\phy$ so as to make
the total number of levels be a multiple of $l$ and we copy
the states at the leaves of $\phy$ to all their descendants in the
new tree. We then apply the estimator as above.

\paragraph{Error Channel.}
We show next that $\bar{s}_\rt$ is a good estimator of the root state
under the conditions of Theorem~\ref{thm:potts}.
Let
\begin{equation*}
\overm{\rt}{l}
= \left(
\prob[\bar{s}_\rt = j\,|\,s_\rt = i]
\right)_{i,j \in \states}.
\end{equation*}
Proposition~\ref{prop:estimpotts} shows
that this ``error channel'' is of the Potts type with
bounded length, no matter how deep the tree.
The key behind our reconstruction algorithm in the next section
will be to think of this error channel as an ``extra edge''
in the Markov model.
\begin{proposition}[Root Estimator from Diluted Trees]\label{prop:estimpotts}
Let $\critperc = \ln 2$.
Then, for all $0 < g < \critperc$ ,
we can find $l > 0$, $\QQ \geq 2$ and $0 < \maxbias < +\infty$ such that
\begin{equation*}
\overm{\rt}{l} = e^{\bias_\rt Q},
\end{equation*}
where $\bias_\rt \leq \maxbias$ and
$Q = Q^{(\qq)}$, for all $\qq > \QQ$
and all $\phy \in \shmgphy^{0,g}$.
\end{proposition}
\begin{proof}
The proof is based on a random cluster argument of Mossel~\cite{Mossel:01}.
Fix $0 < f < g < \critperc$.
In~\cite{Mossel:01}, it is shown that one can choose $\eps > 0$ small enough and $l, R$
large enough
such that
\begin{equation}\label{eq:diluted1}
\prob[\bcal_{i,l}\,|\,s_\rt = i] \geq \eps,
\end{equation}
and
\begin{equation}\label{eq:diluted2}
\prob[\bcal_{i,l}\,|\,s_\rt \neq i] \leq \eps/2,
\end{equation}
for all $\qq > \QQ$
and all $\phy = (V,E,[n],\rt;\weight) \in \shmgphy^{0,g}$.
The proof in~\cite{Mossel:01} actually assumes that all $\weight_e$'s
are equal to $g$. However, the argument
still holds when $\weight_e \leq g$ for all $e$ since smaller $\weight$'s
imply smaller $\delta$'s which can only strenghten inequalities
(\ref{eq:diluted1}) and (\ref{eq:diluted2}) by a standard domination
argument. (For (\ref{eq:diluted2}), see the original argument in~\cite{Mossel:01}.)

Therefore, we have
\begin{eqnarray*}
\overm{\rt}{l}_{ii}
&=& \prob[i \in B_l\,|\,s_\rt = i]\prob[X = i]
+ \frac{1}{\qq - 1}\prob[X \notin B_l\,|\,s_\rt = i, X \neq i]\prob[X \neq i]\\
&\geq& \eps \left(\frac{1}{\qq}\right) + \frac{1}{\qq - 1}(1 - \eps/2) \left(\frac{\qq - 1}{\qq}\right)\\
&=& \frac{1}{\qq} + \frac{\eps}{2\qq}.
\end{eqnarray*}
Also, by symmetry, we have for $i\neq j$
\begin{eqnarray*}
\overm{\rt}{l}_{ij}
&=& \frac{1}{\qq - 1}\left(1 - \overm{\rt}{l}_{ii}\right)\\
&\leq& \frac{1}{\qq} - \frac{\eps}{2\qq(\qq-1)}.
\end{eqnarray*}
Hence, the channel $\overm{\rt}{l}$ is of the form $e^{\bias_\rt Q}$
with $\bias_\rt \leq \maxbias$ where,
by the relation between $\delta$ and $\weight$ given in Example~\ref{ex:symmetric},
we can take
\begin{eqnarray*}
\maxbias
&=& -\ln\left(1 - \qq\left(\frac{1}{\qq} - \frac{\eps}{2\qq(\qq-1)}\right)\right)\\
&=& -\ln\left(\frac{\eps}{2(\qq - 1)}\right).
\end{eqnarray*}
This concludes the proof.
\end{proof}

\subsection{Reconstruction Algorithm}

Our reconstruction algorithm is based on standard distance-based quartet techniques.
Let $\phy = (V,E,[n],\rt;\weight) \in \shmgphy^{f,g}$ be a homogeneous phylogeny
that we seek to reconstruct from $k$ samples of the corresponding Potts
model at the leaves
$(s^i_{[n]})_{i=1}^k \in (\states^{[n]})^k$.

\paragraph{Distances.}
For two nodes $u,v \in V$, we may relate their distance to the probability that their states agree
\begin{equation*}
\weight(u,v) = \sum_{e\in \path(u,v)} \weight_e
= -\ln\left(1 - \left(\frac{\qq}{\qq - 1}\right) \prob[s_u \neq s_v]\right),
\end{equation*}
and so a natural way to estimate $\weight(u,v)$ is to consider
the estimator
\begin{equation*}
\eweight(u,v)
= -\ln\left(1 - \left(\frac{\qq}{\qq - 1}\right)\frac{1}{k}\sum_{i=1}^k\ind\{s^i_u \neq s^i_v\}\right).
\end{equation*}
Of course, given samples at the leaves, this estimator can only be used for $u,v\in [n]$.
Instead, when $u,v$ are internal nodes we first reconstruct their sequence using
Proposition~\ref{prop:estimpotts}. We will then over-estimate the true distance by an amount
not exceeding $2\maxbias$ on average. For $u,v \in V - [n]$, let
\begin{equation*}
\weightb(u,v) = \weight(u,v) + \bias_u + \bias_v,
\end{equation*}
using the notation of Proposition~\ref{prop:estimpotts}.
We also let $\{\bar{s}^i_u\}_{i=1}^k, \{\bar{s}^i_v\}_{i=1}^k$
be the reconstructed states at $u,v$. By convention we let
\begin{equation*}
\weightb(a,b) = \weight(a,b),
\end{equation*}
and
\begin{equation*}
\bar{s}^i_a = s^i_a,\ \forall i=1,\ldots,k,
\end{equation*}
for $a,b \in [n]$.
Note that, at the beginning of the algorithm,
the phylogeny is not known,
making it impossible to compute $\{\bar{s}^i_u\}_{i=1}^k$ for internal
nodes.
However as we reconstruct parts of the tree we will progressively
compute the estimated sequences of uncovered internal nodes.

By standard concentration inequalities, $\weightb(u,v)$ can be well approximated
with $k = O(\log n)$ as long as $\weightb(u,v) = O(1)$.
For $u,v\in V$ let
\begin{equation*}
\eweight(u,v)
= -\ln\left(1 - \left(\frac{\qq}{\qq - 1}\right)\frac{1}{k}\sum_{i=1}^k\ind\{\bar{s}^i_u \neq \bar{s}^i_v\}\right).
\end{equation*}
Recall the notation of Example~\ref{ex:homo}.
\begin{lemma}[Distorted Metric: Short Distances~\cite{ErStSzWa:99a}]\label{lem:distmet1}
Let $0\leq h'< h$ and let $u,v \in L^{(h)}_{h'}$ be distinct leaves.
For all $D > 0$, $\delta > 0$, $\gamma > 0$,
there exists $c =
c(D, \delta,\gamma) > 0$, such that if the following conditions hold:
\begin{itemize}
\item $\mathrm{[Small\ Diameter]}$ $\weightb(u,v) < D$,
\item $\mathrm{[Sequence\ Length]}$ $k = c' \log n$ for $c' > c$,
\end{itemize}
then
\begin{equation*}
\left|\weightb(u,v)-\eweight(u,v)\right|< \delta,
\end{equation*}
with probability at least $1-n^{-\gamma}$.
\end{lemma}
\begin{lemma}[Distorted Metric: Diameter Test~\cite{ErStSzWa:99a}]\label{lem:distmet2}
Let $0\leq h'< h$ and $u,v \in L^{(h)}_{h'}$.
For all $D > 0$, $W > 5$, $\gamma > 0$,
there exists $c =
c(D, W,\gamma) > 0$, such that if the following conditions hold:
\begin{itemize}
\item $\mathrm{[Large\ Diameter]}$ $\weightb(u,v) > D + \ln W$,
\item $\mathrm{[Sequence\ Length]}$ $k = c' \log n$ for $c' > c$,
\end{itemize}
then
\begin{equation*}
\eweight(u,v)
> D + \ln \frac{W}{2},
\end{equation*}
with probability at least $1-n^{-\gamma}$.
On the other hand, if the first condition above is replaced by
\begin{itemize}
\item $\mathrm{[Small\ Diameter]}$ $\weightb(u,v)
< D + \ln \frac{W}{5}$,
\end{itemize}
then
\begin{equation*}
\eweight(u,v)
\leq D + \ln \frac{W}{4},
\end{equation*}
with probability at least $1-n^{-\gamma}$.
\end{lemma}

\paragraph{Quartest Tests.}
Let $0 \leq h' < h$ and
$\qcal_0 = \{a_0,b_0,c_0,d_0\} \subseteq L^{(h)}_{h'}$.
The topology of $T^{(h)}$ restricted to $\qcal_0$
is completely characterized by a bi-partition or
{\em quartet split} $q_0$ of the form:
$a_0 b_0 | c_0 d_0$, $a_0 c_0 | b_0 d_0$ or
$a_0 d_0 | b_0 c_0$.
The most basic operation in quartet-based reconstruction
algorithms is the inference of such quartet splits.
In distance-based methods in particular, this is usually done
by performing the so-called {\em four-point test}:
letting
\begin{equation*}
\fcal(a_0 b_0 | c_0 d_0)
= \frac{1}{2}[\weight(a_0,c_0) + \weight(b_0,d_0)
- \weight(a_0,b_0) - \weight(c_0,d_0)],
\end{equation*}
we have
\begin{equation*}
q_0
=
\left\{
\begin{array}{ll}
a_0 b_0 | c_0 d_0 & \mathrm{if\ }\fcal(a_0,b_0|c_0,d_0) > 0\\
a_0 c_0 | b_0 d_0 & \mathrm{if\ }\fcal(a_0,b_0|c_0,d_0) < 0\\
a_0 d_0 | b_0 c_0 & \mathrm{o.w.}
\end{array}
\right.
\end{equation*}
Note that adding ``extra edges'' at the nodes $a_0, b_0, c_0, d_0$ as implied
in Proposition~\ref{prop:estimpotts} does not affect the topology
of the quartet.

Since Lemma~\ref{lem:distmet1} applies only to short distances, we also
perform a diameter test. We let $\widehat\fcal(a_0 b_0 | c_0 d_0) = +\infty$
if $\max_{u,v \in \qcal_0} \eweight(u,v) > D + \ln \frac{W}{4}$
and otherwise
\begin{equation*}
\widehat\fcal(a_0 b_0 | c_0 d_0)
=
\frac{1}{2}[\eweight(a_0,c_0) + \eweight(b_0,d_0)
- \eweight(a_0,b_0) - \eweight(c_0,d_0)].
\end{equation*}
Finally we let
\begin{equation*}
\deep(a_0,b_0|c_0,d_0) = \ind\{\widehat\fcal(a_0 b_0 | c_0 d_0) > f/2\}.
\end{equation*}

\noindent\textbf{Algorithm.}
The algorithm is detailed in Figure~\ref{fig:algo}.
The proof of its correctness is left to the reader.
This concludes the proof of Theorem~\ref{thm:potts}.

\begin{figure*}[!ht]
\framebox{
\begin{minipage}{12.2cm}
{\small \textbf{Algorithm}\\
\textit{Input:} Sequences $(s^i_{[n]})_{i=1}^k \in (\states^{[n]})^k$;\\
\textit{Output:} Tree;

\begin{itemize}
\item Let $\zcal_{0}$ be the set of leaves.
\item For $h' = 0,\ldots,h-1$,
\begin{enumerate}
\item \textbf{Four-Point Test.}
Let
\begin{equation*}
\rcal_{h'} = \{q = ab|cd\ :\ \forall a,b,c,d \in \zcal_{h'}\ \text{distinct such that}\ \deep(q) = 1\}.
\end{equation*}

\item \textbf{Cherries.} Identify the cherries in $\rcal_{h'}$, that is, those pairs of vertices
that only appear on the same side of the quartet splits in $\rcal_{h'}$.
Let
\begin{equation*}
\zcal_{h'+1} = \{a_1^{(h'+1)},\ldots,a_{2^{h - (h'+1)}}^{(h'+1)}\},
\end{equation*}
be the parents of the cherries in $\zcal_{h'}$.

\item \textbf{Reconstructed Sequences.}
For all $u \in \zcal_{h'+1}$, compute
$(\bar{s}^i_u)_{i=1}^k$.

\end{enumerate}

\end{itemize}

}
\end{minipage}
} \caption{Algorithm.} \label{fig:algo}
\end{figure*}

\section{General Lower Bound}\label{section:lower}

Here we prove the following statement which implies Theorem~\ref{thm:lower}:
\begin{theorem}[Polynomial Lower Bound on PTR] \label{thm:coupling}
Consider the phylogenetic reconstruction problem for homogeneous trees with fixed edge length $\weight(e) = \weight > 0$ for all edges $e\in E$.
Assume further that the ASR problem for edge length $\weight$ and matrix $Q$ is not solvable and that moreover
$\weight > \critks$. Then there exists $\alpha = \alpha(\weight) > 0$ such that the probability of correctly reconstructing
the tree is at most $O(n^{-\alpha})$ assuming $k \leq n^\alpha$.
\end{theorem}

For general mutation rates $Q$, it is not known if there is a {\em unique} reconstruction threshold $\critmlq$ such that
ASR is possible for $\weight < \critmlq$ and impossible for $\weight > \critmlq$.
For models for which such a threshold exists Theorem~\ref{thm:coupling} above shows the
impossibility of phylogenetic reconstruction for $\weight > \critmlq$.
The existence of the threshold $\critmlq$ has been established for a few models,
e.g.~for so-called random cluster models,
which include the binary asymmetric channel and the Potts model~\cite{Mossel:01}.

The proof of Theorem~\ref{thm:coupling} is based on the following two lemmas.
It is useful to write $n = 2^{\ell}$ for the number of leaves of a homogeneous
tree with $\ell$ levels.
\begin{lemma}[Reconstructing a Deep Subtree] \label{lem:coup1}
Consider the PTR problem for homogeneous trees with fixed edge length $\weight$.
Let $\mu_Q^{\ell,i}$ denote the distribution at the leaves on a homogeneous $\ell$-level tree
with fixed edge length $\weight$, root value $i$, and rate matrix $Q$.
Suppose there exists a number $0 < \alpha < 1$ such that for every $\ell$ and all $i$ one can write
$\mu_Q^{\ell,i} = (1-\eps) \bar{\mu} + \eps \mu'^{i}$
for some probability measures $\mu'^{{i}}, i \in \states$, $\bar{\mu}$,
and $\eps = O(2^{-\alpha \ell})$.
Then the probability of correctly reconstructing homogeneous phylogenetic trees
with edge length $\weight$ assuming
$k \leq n^{\alpha/10}$ is at most $O(n^{-\alpha/2})$.
\end{lemma}
\begin{lemma}[Leaf Distribution Decomposition] \label{lem:coup2}
Consider the ASR problem for homogeneous trees with fixed edge length $\weight$.
Assume further that the ASR problem for $Q$ with edge length $\weight$ is not solvable and
further $\weight > \critks$.
Then there is an $\alpha = \alpha(\weight) > 0$ for which the following holds.
There exists a sequence $\eps_{\ell} = O(2^{-\alpha \ell})$ such that for all $i\in \states$ one can write
$\mu_Q^{\ell,i} = (1-\eps) \bar{\mu} + \eps \mu'^{i}$ for some probability measures $\mu'^{i}, i \in \states$ and $\bar{\mu}$.
\end{lemma}
\begin{app-proof}{Lemma~\ref{lem:coup1}}
Let $r$ be chosen so that $2^{r-1} < n^{\alpha/20} \leq 2^{r}$. (Note that $r < \ell$.)
Consider the following distribution: first, pick a homogeneous tree $T$ on $\ell$ levels, where the first $r$
levels are chosen uniformly at random among $r$-level homogeneous trees and the remaining levels are fixed (i.e., deterministic);
second, pick $k$ samples of a Markov model with rate matrix $Q$ and fixed edge length $\weight$ on the resulting tree.

Let $\acal$ be a phylogenetic reconstruction algorithm.
Our goal is to bound the success probability of $\acal$ on the random model above.
We may assume that the bottom $\ell - r$ levels are given to $\acal$ (as it may ignore this information)
and that $\acal$ is deterministic
(as a simple convexity argument shows that deterministic algorithms achieve the highest success probability).

Note that the assumption of the lemma implies that, for a single sample,
we can {\em simultaneously} couple the distribution
at the leaves of all the given subtrees of $\ell-r$ levels---except with probability $O(2^r 2^{-\alpha(\ell - r)}) = O(n^{-9\alpha/10})$.
This can be achieved by starting the coupling at level $r$ (from the root) of the tree.
Repeating this for the $n^{\alpha/10}$ samples we obtain the following.
Let $\mu_T$ denote the measure on the $n^{\alpha/10}$ samples
at leaves of $T$.
Then there exists measures $\mu, \mu'_T$ and $\eps = O(n^{-8 \alpha/10})$ such that
$\mu_T = (1-\eps) \mu + \eps \mu'_T$.

Write $N_{r}$ for the number of leaf-labelled complete binary trees on $r$ levels. Write $\ecal(s,\acal,T)$ for the indicator of the event that
the $k$ samples are given by $s$ and that $\acal$ recovers $T$.
The success probability of $\acal$ is then given by
{\small
\begin{eqnarray}
&&\sum_T N_r^{-1}  \left(\sum_{s} \mu_T(\ecal(s,\acal,T))\right) \nonumber\\
&& \qquad =
(1-\eps) N_r^{-1} \sum_{s} \sum_T \mu(\ecal(s,\acal,T))
+ \eps N_r^{-1} \sum_T \sum_{s} \mu'_T(\ecal(s,\acal,T)).\label{eq:success}
\end{eqnarray}
}
For the second term note that
\[
\sum_{s} \mu'_T(\ecal(s,\acal,T)) \leq \sum_{s} \mu'_T(s) = 1,
\]
and therefore the second term in (\ref{eq:success}) is bounded by $\eps$.
Furthermore for each $s$, $\sum_T \mu(\ecal(s,\acal,T)) = \mu(s)$
by definition
and $\sum_{s} \mu(s) = 1$ so the first term in (\ref{eq:success}) is bounded by $(1-\eps) N_r^{-1}$.

Thus overall, the bound on the probability of correct reconstruction is $\eps + (1-\eps) N_r^{-1}$. Using the facts that
$N_r = \Omega(2^{2^r}) = \Omega(2^{n^{0.1\alpha}}) = \Omega(n^{\alpha/2})$
and $\eps = O(n^{-8\alpha/10})$ concludes the proof.
\end{app-proof}
\begin{app-proof}{Lemma~\ref{lem:coup2}}
For $\delta > 0$ and $r' > 0$, let $\mu_Q^{\ell-r',i}(\delta)$ be the same measure as $\mu_Q^{\ell-r',i}$,
except that, for each leaf,
independently with probability $1-\delta$, the state at the leaf is replaced by $*$ (which does not belong
to the original alphabet).
The key to the proof is the main result of~\cite{JansonMossel:04} where it is shown that if $\weight > \critks$
then the following holds:
There exist fixed $\delta > 0, \alpha > 0$ such that
\begin{equation} \label{eq:janson}
\mu_Q^{\ell-r',i}(\delta) = (1-\eps) \bar{\mu}(\delta) + \eps \mu'^{i}(\delta),
\end{equation}
where $\eps = O(2^{-\alpha (\ell-r')})$
for some probability measures $\mu'^{i}(\delta)$ and $\bar{\mu}(\delta)$.

The fact that there is no reconstruction (ASR)  at edge length $\weight$ implies that there exists a fixed $r'$ and measures
$\bar{\nu}$ and $\nu'^i$ such that
\[
\mu_Q^{r',i} = (1-\delta) \bar{\nu} + \delta \nu'^i.
\]
This implies in particular that we can simulate the mutation process on an $\ell$-level tree by first using
the measure
$\mu_Q^{\ell,i}(\delta)$ and then applying the following rule: for each node $v$ at level $\ell-r'$ independently
\begin{itemize}
\item
If the label at $v$ is $*$ then generate the leaf states on the subtree rooted at $v$ according to the measure $\bar{\nu}$.
\item
Else if it is labeled by $i$, sample leaf states on the subtree below $v$ from the measure $\nu'^i$.
\end{itemize}
The desired property of the measures $\mu_Q^{\ell,i}$ now follows
from the fact that the measures $\mu_Q^{\ell,i}(\delta)$ have the desired property by~(\ref{eq:janson}).
\end{app-proof}

%

\bibliographystyle{alpha}
\bibliography{thesis}

\end{document}